\documentclass[12pt]{amsart}
\usepackage{amssymb,amsthm,amsmath,epsfig,latexsym}
\usepackage{calc}
\usepackage{graphicx}
\usepackage{psfrag}

\begin{document}

\newcommand{\mmbox}[1]{\mbox{${#1}$}}
\newcommand{\proj}[1]{\mmbox{{\mathbb P}^{#1}}}
\newcommand{\Cr}{C^r(\Delta)}
\newcommand{\CR}{C^r(\hat\Delta)}
\newcommand{\affine}[1]{\mmbox{{\mathbb A}^{#1}}}
\newcommand{\Ann}[1]{\mmbox{{\rm Ann}({#1})}}
\newcommand{\caps}[3]{\mmbox{{#1}_{#2} \cap \ldots \cap {#1}_{#3}}}
\newcommand{\N}{{\mathbb N}}
\newcommand{\Z}{{\mathbb Z}}
\newcommand{\R}{{\mathbb R}}
\newcommand{\Tor}{\mathop{\rm Tor}\nolimits}
\newcommand{\Ext}{\mathop{\rm Ext}\nolimits}
\newcommand{\Hom}{\mathop{\rm Hom}\nolimits}
\newcommand{\im}{\mathop{\rm Im}\nolimits}
\newcommand{\rank}{\mathop{\rm rank}\nolimits}
\newcommand{\supp}{\mathop{\rm supp}\nolimits}
\newcommand{\arrow}[1]{\stackrel{#1}{\longrightarrow}}
\newcommand{\CB}{Cayley-Bacharach}
\newcommand{\coker}{\mathop{\rm coker}\nolimits}
\sloppy
\newtheorem{defn0}{Definition}[section]
\newtheorem{prop0}[defn0]{Proposition}
\newtheorem{conj0}[defn0]{Conjecture}
\newtheorem{thm0}[defn0]{Theorem}
\newtheorem{lem0}[defn0]{Lemma}
\newtheorem{corollary0}[defn0]{Corollary}
\newtheorem{example0}[defn0]{Example}

\newenvironment{defn}{\begin{defn0}}{\end{defn0}}
\newenvironment{prop}{\begin{prop0}}{\end{prop0}}
\newenvironment{conj}{\begin{conj0}}{\end{conj0}}
\newenvironment{thm}{\begin{thm0}}{\end{thm0}}
\newenvironment{lem}{\begin{lem0}}{\end{lem0}}
\newenvironment{cor}{\begin{corollary0}}{\end{corollary0}}
\newenvironment{exm}{\begin{example0}\rm}{\end{example0}}

\newcommand{\defref}[1]{Definition~\ref{#1}}
\newcommand{\propref}[1]{Proposition~\ref{#1}}
\newcommand{\thmref}[1]{Theorem~\ref{#1}}
\newcommand{\lemref}[1]{Lemma~\ref{#1}}
\newcommand{\corref}[1]{Corollary~\ref{#1}}
\newcommand{\exref}[1]{Example~\ref{#1}}
\newcommand{\secref}[1]{Section~\ref{#1}}

\newcommand{\std}{Gr\"{o}bner}
\newcommand{\jq}{J_{Q}}



\title {Topological Criteria for $k-$Formal Arrangements}
\author[\c Stefan O. Toh\v aneanu]{\c Stefan O. Toh\v aneanu$^1$}
\address{Department of Mathematics,
Texas A\&M University, College Station, TX 77843}
\email{tohanean@math.tamu.edu}

\thanks{$^1$Partially supported by NSF DMS 03-11142 and ATP 010366-0103.}


\begin{abstract}
\noindent We prove a criterion for $k-$formality of arrangements,
using a complex constructed from vector spaces introduced in
\cite{bt}. As an application, we give a simple description of
$k-$formality of graphic arrangements: Let $G$ be a connected
graph with no loops or multiple edges. Let $\Delta$ be the flag
(clique) complex of $G$ and let $H_{\bullet}(\Delta)$ be the
homology of the chain complex of $\Delta$. If $\mathcal A_G$ is
the graphic arrangement associated to $G$, we will show that
$\mathcal A_G$ is $k-$formal if and only if $H_i(\Delta)=0$ for
every $i=1,\ldots,k-1$.
\end{abstract}
\maketitle

\section{Introduction}

In \cite{fr}, Falk and Randell introduced the notion of a formal arrangement. An arrangement is formal iff every
linear dependency among the defining forms of the hyperplanes can be expressed as linear combination of
dependencies among exactly 3 defining forms. Many interesting classes of arrangements are formal: in \cite{fr},
Falk and Randell proved that $K(\pi,1)$ arrangements and arrangements with quadratic Orlik-Solomon ideal are
formal and, in \cite{y}, Yuzvinsky showed that free arrangements are also formal; and gave an example showing
that formality does not depend on the intersection lattice. In \cite{bt}, Brandt and Terao generalized the
notion of formality to $k-$formality, proving that every free arrangement is $k-$formal. For this they
introduced the concept of 'higher' relation spaces, which capture 'the dependencies among dependencies'.

In the first section of this paper we briefly recall the notions
of relation spaces and $k-$formality. By rewriting the
definitions, we obtain a lemma characterizing $k-$formality
topologically. Then we apply this criteria for graphic
arrangements to obtain a description of $k-$formality in terms of
the homology of a certain chain complex. With this it is easy to
produce examples of graphic arrangements which are $k-$formal but
not $(k+1)-$formal, for any given $k$.

\section{Preliminaries}

In what follows we adopt all the notation from \cite{bt}. Let
$\mathcal A$ be an arrangement of $n$ hyperplanes in a vector
space $V$ over a field $\mathbb K$. For each $H\in \mathcal A$ we
fix the defining form $\alpha_H\in V^*$.

Define a map $\phi: E(\mathcal A):=\oplus_{H\in \mathcal A}\mathbb
K e_H \rightarrow V^*$, by $\phi(e_H)=\alpha_H$, where $E(\mathcal
A)$ is the vector space with basis $\{e_H\}$.

Let $F(\mathcal A)$ be the kernel of this map. Then $\dim
F(\mathcal A)=n-r(\mathcal A)$ where $r(\mathcal A)$ is the rank
of $\mathcal A$. The vector space $F(\mathcal A)$ describes which
linear forms are linearly dependent, as well as the dependency
coefficients (up to scalar multiplication). We will refer to
elements of $F(\mathcal A)$ as {\em relations}.

Let $F_2(\mathcal A)$ be the subspace of $F(\mathcal A)$ generated
by the relations corresponding to dependencies of exactly 3 linear
forms.

\begin{defn} $\mathcal A$ is formal iff $F(\mathcal
A)=F_2(\mathcal A)$.
\end{defn}

\begin{defn} For $3\leq k \leq r(\mathcal A)$, recursively define
$R_k(\mathcal A)$ to be the kernel of the map $$\pi_{k-1} =
\pi_{k-1}(\mathcal A): \bigoplus_{X\in L, r(X)=k-1}R_{k-1}
(\mathcal A_X) \rightarrow R_{k-1}(\mathcal A), $$ where $L$ is
the lattice of intersections of $\mathcal A$ and $\pi_{k-1}$ is
the sum of the inclusion maps $R_{k-1}(\mathcal
A_X)\hookrightarrow R_{k-1}(\mathcal A)$. We identify
$R_2(\mathcal A)$ with $F(\mathcal A)$.
\end{defn}

To simplify notation, for $k\geq 2$ we will denote with
$D_k=D_k(\mathcal A)$ the vector space $\bigoplus_{X\in L,
r(X)=k}R_{k} (\mathcal A_X)$.

\begin{defn} We define

\begin{enumerate}
    \item An arrangement is 2-formal if it is formal.
    \item For $k\geq 3$, $\mathcal A$ is $k-$formal iff it is
    $(k-1)-$formal and the map $\pi_k:D_k\rightarrow R_k(\mathcal
    A)$ is surjective.
\end{enumerate}
\end{defn}

\vskip .2in

\begin{lem} For any arrangement $\mathcal A$, the following sequence
of vector spaces and maps form a complex:
$$D_{\bullet}: 0\longrightarrow \cdots \stackrel{d_3}\longrightarrow
D_2 \stackrel{d_2} \longrightarrow D_1 \stackrel{d_1} \longrightarrow D_0 \longrightarrow 0,$$ where $D_0=V^*$,
$D_1=E(\mathcal A)$ and for $k\geq 2$, $D_k$ are the spaces from the notations above. Also, $d_1=\phi$ and
$d_k:D_k\rightarrow D_{k-1}, d_k=\pi_k$ for $k\geq 2$.
\end{lem}
\begin{proof}
We have $d_k(D_k)=\pi_k(D_k)\subseteq R_k(\mathcal
A)=ker(\pi_{k-1})\subseteq D_{k-1}$. So $d_k$ is well defined.
Also, $d_{k-1}\circ d_k(v)= \pi_{k-1}(\pi_k(v))=0$ for any $v\in
D_k$, as $\pi_k(v)\in R_k(\mathcal A)=ker(\pi_{k-1})$. So, indeed
we have a complex.
\end{proof}
\begin{lem} $\mathcal A$ is $k-$formal iff $H_i(D_{\bullet})=0$ for every
$i=1,\ldots,k-1$
\end{lem}
\begin{proof} $\pi_l$ is surjective iff $\forall w\in R_l(\mathcal
A)$ there exists $v\in D_l$ such that $\pi_l(v)=w$.

We have $R_l(\mathcal A)=ker(\pi_{l-1})=ker(d_{l-1})$ and $w=\pi_l(v)=d_l(v)\in Im(d_l)$. So we get
$ker(d_{l-1})\subseteq Im(d_l)$ which give us $H_{l-1}(D)=0$.
\end{proof}

\begin{exm} In this example we will discuss \cite{bt}, Example
5.1., in terms of the homology of the above complex. We must specify that all the computations are already done
in \cite{bt}, and we are just translating into topological language.

$\mathcal A$ is a real essential arrangement of rank 4 consisting of 10 hyperplanes, defined by the vanishing of
the following linear forms: $\alpha_1=x_3,\alpha_2=x_3-x_4,\alpha_3=x_2,\alpha_4=x_2+x_3-
2x_4,\alpha_5=x_1,\alpha_6=x_1+x_3-2x_4,\alpha_7=x_2+2x_3-2x_4,
\alpha_8=x_1+2x_3-2x_4,\alpha_9=x_1+x_2+x_3-2x_4,\alpha_{10}=x_4$.

So $D_0=\mathbb R^4$, $D_1=\mathbb R^{10}$ and the map
$d_1:D_1\longrightarrow D_0$ is just the map $\phi$ and has rank
4. Therefore $ker(d_1)$ has dimension $10-4=6$.

We have 7 nondegenerate rank 2 elements in $L(\mathcal A)$ and
each is an intersection of exactly 3 hyperplanes. So we have 7
relations of length 3: $\alpha_1-\alpha_2-\alpha_{10}=0, \alpha_1
+\alpha_4-\alpha_7=0, \alpha_1+\alpha_6-\alpha_8=0, 2\alpha_2+
\alpha_3-\alpha_7=0, 2\alpha_2+\alpha_5-\alpha_8=0, \alpha_3+
\alpha_6-\alpha_9=0, \alpha_4+\alpha_5-\alpha_9=0$.

Therefore $D_2=\mathbb R^7$. The matrix of the map
$d_2:D_2\longrightarrow D_1$ is exactly the matrix in \cite{bt},
page 61 \begin{displaymath} \left(
\begin{array}{cccccccccc}
1&-1&0&0&0&0&0&0&0&-1\\
1&0&0&1&0&0&-1&0&0&0\\
1&0&0&0&0&1&0&-1&0&0\\
0&2&1&0&0&0&-1&0&0&0\\
0&2&0&0&1&0&0&-1&0&0\\
0&0&1&0&0&1&0&0&-1&0\\
0&0&0&1&1&0&0&0&-1&0
\end{array}
\right),
\end{displaymath} and it has rank 6. So $\dim Im(d_2)=6$ and $\dim
ker(d_2)=7-6=1$.

Also in \cite{bt} we have listed all the elements of rank 3 from
$L(\mathcal A)$: $\{1,2,9,10\}$, $\{3,6,9,10\}$, $\{4,5,9,10\}$,
$\{1,3,6,8,9\}$, $\{1,4,5,7,9\}$, $\{1,4,6,7,8\}$,
$\{2,3,5,7,8\}$, $\{2,3,6,7,9\}$, $\{2,4,5,8,9\}$,
$\{3,4,5,6,9\}$, $\{1,2,3,4,7,10\}$, $\{1,2,5,6,8,10\}$.

If $X$ is such an element (with $r(X)=3$), then $R_3(\mathcal A_X)\neq 0$ means that there is at least a
relation among the relations of length 3 of elements of rank 2 in $L(\mathcal A_X)$. The nondegenerate rank 2
elements in $L(\mathcal A_X)$ are nondegenerate rank 2 elements in $L(\mathcal A)$ and these are listed above.
It is not difficult to check which are the relations of length 3 for each rank 3 element in $\mathcal A$. For
reference, these are listed in the chart on page 62 in \cite{bt}. Also, there is no problem to check that for
each $r(X)=3$, the length 3 relations are linearly independent. Therefore we conclude that $D_3=0$.

So the complex we get is:
$$D_{\bullet}: 0\longrightarrow \mathbb R^7 \longrightarrow \mathbb R^{10}
\longrightarrow \mathbb R^4 \longrightarrow 0$$ with homology:
$H_1(D_{\bullet})=0$ and $H_2(D_{\bullet})=1$. So $\mathcal A$ is
formal, but not 3-formal.
\end{exm}

\vskip .2in

Graphic arrangements are a class of arrangements possessing many nice properties (see, for example, \cite{KS06},
\cite{SS02}, \cite{SS06}, \cite{St72}), and for this class there is a pleasant combinatorial interpretation of
Lema 2.5. In the next section, $G$ denotes a connected graph with no loops or multiple edges. For the graphic
arrangement $\mathcal A_G$, we will identify the complex above with the chain complex of the flag complex of
$G$. Then, with Lemma 2.5., the statement in the abstract will be natural.

\section{Graphic Arrangements}

Let $G$ be a connected graph on vertices $[n]=\{1,\ldots,n\}$ with
no loops or multiple edges. The \textit{flag(clique) complex}
$\Delta=\Delta(G)$ is the simplicial complex with:
\begin{itemize}
    \item The 0-faces = the vertices of $G$.
    \item The 1-faces = the edges of $G$.
    \item For $i\geq 2$, the $i-$faces = the $K_{i+1}$ (i.e., complete
    graph on $i+1$ vertices) subgraphs of $G$.
\end{itemize}

For $i\geq 0$, let $a_i$ be the number of $i-$faces of $\Delta$.
We have the natural chain complex of $\Delta$: $$ 0\longrightarrow
\cdots \stackrel{f_3}\longrightarrow C_2 \stackrel{f_2}
\longrightarrow C_1 \stackrel{f_1} \longrightarrow C_0
\longrightarrow 0,$$ where $C_i=\mathbb K^{a_i}$ and $f_i:C_i
\rightarrow C_{i-1}$ is the usual differential given in terms of
generators: $f_i([n_1,\ldots,n_{i+1}]) = \sum_{j=1}^{i+1}
(-1)^{j-1}[n_1,\ldots,\hat{n_j},\ldots,n_{i+1}]$.

The homology of this complex will be denoted by
$H_{\bullet}(\Delta)$.

\vskip .2in

By definition, the graphic arrangement associated to $G$ is
$\mathcal A = \mathcal A_G =\{ker\{\alpha_{ij}\}
|\alpha_{ij}=x_i-x_j, i<j \mbox{ and } [ij] \mbox{ is an edge in }
G\}$. Note that $\mathcal A$ is an arrangement in $V=\mathbb
K^{a_0}$ of rank $a_0-1$ (if $G$ is connected) and consists of
$a_1$ (= the number of edges in $G$) hyperplanes.

Notice that from the beginning we fixed the defining forms
$\alpha_{ij}$. To be consistent with notation, $e_{ij}, i<j$ will
be the symbols in $E(\mathcal A)$ (i.e., $\phi(e_{ij})
=\alpha_{ij}$). With these, we can identify $D_1=E(\mathcal A)$
with $C_1$ by $e_{ij}\leftrightarrow [ij]$ for $i<j$.

\vskip .2in

If we fix the form of the elements in the basis of $D_i$'s and with proper notations of those, the
correspondence between the two complexes will become natural. The next lemma will do this, but before we state
and prove it here is the flavor of it:

For $X\in L$, let $G_X$ be the subgraph of $G$ built on the edges
corresponding to the hyperplanes in $X$.

We have $D_2=\oplus_{X\in L_2}R_2(\mathcal A_X)$. Suppose for an $X\in L_2$ we have $R_2(\mathcal
A_X)=F(\mathcal A_X)\neq 0$. This means that we must have a dependency (relation) among some of the linear forms
corresponding to some edges in $G_X$. But this translates in the fact that $G_X$ contains a cycle. If the length
of this cycle is $\geq 4$, then the linear forms corresponding to 3 consecutive edges in the cycle are linearly
independent. This contradicts the fact that $rk(X)=2$. So $G_X$ contains a triangle. If we have an extra edge in
$G_X$, beside those from the triangle, then  the linear form of this extra edge and the linear forms associated
to two of the edges of the triangle are linearly independent. Again we get a contradiction with the fact that
$rk(X)=2$. So $G_X=$ a triangle. So each nonzero summand of $D_2$ corresponds to a triangle in $G$. The converse
of this statement is obvious.

\begin{lem}(The Recursive Identification Lemma) Let $X\in L$ with $r(X)=l,
l\geq 2$. Then $R_l(\mathcal A_X)\neq 0$ iff $G_X$ is a $K_{l+1}$
subgraph of $G$. More, $\dim R_l(\mathcal A_X)=1$ and if $G_X=
[i_1i_2\cdots i_{l+1}], i_1<i_2<\cdots i_{l+1}$, then we can pick
a 'special' basis element of $R_l(\mathcal A_X)$ to be the
relation on the special elements corresponding to the $K_l$
subgraphs of $G_X$: $r_{i_2\cdots i_{l+1}}-r_{i_1i_3\cdots
i_{l+1}}+\cdots +(-1)^lr_{i_1i_2\cdots i_l}$. This element is
denoted with $r_{i_1i_2\cdots i_{l+1}}$.
\end{lem}
\begin{proof} Suppose $R_l(\mathcal A_X)\neq 0$. We will use induction on $l$.

For $l=2$ we already seen this case above.

Suppose $l\geq 3$.

By definition, we have $R_l(\mathcal A_X)=ker(\pi_{l-1})$, where
$$ \pi_{l-1}:D_{l-1}(\mathcal A_X)= \bigoplus_{Y\in L(\mathcal A_X),
r(Y)=l-1}R_{l-1}((\mathcal A_X)_Y)\longrightarrow R_{l-1}(\mathcal
A_X).$$

The induction hypothesis is telling that for each $Y\in L(\mathcal
A_X), r(Y)=l-1$ such that $R_{l-1}((\mathcal A_X)_Y)\neq 0$, $G_Y
= [i_1i_2\cdots i_l]$ is a $K_l$ subgraph of $G_X$ and $\dim
R_{l-1} ((\mathcal A_X)_Y)=1$ with $r_{i_1i_2\cdots i_l}$
'special' basis element of $R_{l-1}((\mathcal A_X)_Y)$.

>From this we get first that $\dim D_{l-1}(\mathcal A_X)=$ the
number of $K_l$ subgraphs of $G_X$.

\vskip .2in

The condition $R_l(\mathcal A_X)\neq 0$ is telling us that, since
$R_l(\mathcal A_X)\subseteq D_{l-1}(\mathcal A_X)$, $G_X$ has at
least one $K_l$ subgraph.

If $G_X$ has just one $K_l$ subgraph, then $R_l(\mathcal A_X)=
D_{l-1}(\mathcal A_X)$. But $\pi_{l-1}$ is a sum of inclusions,
and in this particular case it will be exactly an inclusion. So we
get that $R_l(\mathcal A_X)=ker(\pi_{l-1})=0$, which is a
contradiction.

Therefore, $G_X$ has at least two $K_l$ subgraphs.

Let's take two of them $K_l^1$ and $K_l^2$, and first suppose they
do not share any vertex. Let $v\in K_l^1$ and $w\in K_l^2$ be two
vertices of $G_X$. Through $v$ pass exactly $l-1$ edges and the
corresponding linear forms $\alpha_1,\ldots,\alpha_{l-1}$ are
linearly independent. Let's take two edges $[w,w_1]$ and $[w,w_2]$
of $K_l^2$ and let $\beta_1$ and $\beta_2$ be the corresponding
linear forms. Then, $\alpha_1,\ldots,\alpha_{l-1},\beta_1,
\beta_2$ are linearly dependent if at least one of the vertices
$\{w,w_1,w_2\}$ is a vertex in $K_l^1$. Contradiction. Therefore,
$\alpha_1,\ldots,\alpha_{l-1},\beta_1, \beta_2$ are linearly
independent. But this will contradict $r(\mathcal A_X)=l$.

Hence, $K_l^1$ and $K_l^2$ have at least a common vertex $v$.
Suppose $w_1,w_2$ are two vertices of $K_l^2$ but not of $K_l^1$.
Then, through $v$ pass at least $l+1$ edges: $l-1$ from $K_l^1$
and $[v,w_1]$, $[v,w_2]$ from $K_l^2$. The corresponding linear
forms are linearly independent and again we obtain a contradiction
with the fact that $r(\mathcal A_X)=l$.

The conclusion of all of above is that any two distinct $K_l$
subgraphs of $G_X$ have exactly $l-1$ vertices in common. $(*)$

\vskip .2in

Suppose $G_X$ has exactly two $K_l$ subgraphs:
$[1,2,\ldots,l-1,l]$ and $[1,2,\ldots,l-1,l+1]$. Let $r\in
R_l(\mathcal A_X), r\neq 0$. Then $r=r_{1,2,\ldots ,l-1,l}+b
r_{1,2,\ldots,l-1,l+1}$ for some $b\in \mathbb K-\{0\}$. We have
$\pi_{l-1}(r)=0$ in $D_{l-2}(\mathcal A_X)$.  So we get a relation
on the 'special' basis elements of $D_{l-2}(\mathcal A_X)$: $$
0=(r_{2,\ldots ,l-1,l}-r_{1,3,\ldots ,l-1,l}+\cdots
+(-1)^{l-1}r_{1,2,\ldots ,l-1})$$ $$+b(r_{2,\ldots,l-1,l+1}-
r_{1,3,\ldots,l-1,l+1}+\cdots +(-1)^{l-1}r_{1,2,\ldots,l-1}).$$
Observe that this equation is impossible.

So $G_X$ has at least three distinct $K_l$ subgraphs:
$K_l^1=[1,2,\ldots,l-1,l]$, $K_l^2=[1,2,\ldots,l-1,l+1]$ and
$K_l^3$. If both $l$ and $l+1$ are vertices in $K_l^3$, then $l$
and $l+1$ are connected in $G_X$, so $G_X$ contains a $K_{l+1}$
subgraph. If, for example, $l\notin K_l^3$, then from $(*)$ and
since $K_l^i,i=1,2,3$ are distinct we get that
$K_l^3=[1,2,\ldots,l-1,l+2]$, for some other vertex $l+2$ in
$G_X$. Observe that through the vertex $1$ pass at least $l+1$
edges of $G_X$: $[1,2],[1,3],\ldots,[1,l-1],[1,l],[1,l+1],
[1,l+2]$. The corresponding linear forms of these edges are
linearly independent so we get a contradiction with the fact that
$r(\mathcal A_X)=l$.

We can conclude that $G_X$ contains a $K_{l+1}$ subgraph. Now, if
there exists an extra edge of $G_X$ not on this $K_{l+1}$, then
the corresponding linear form of this edge together with the
corresponding linear forms of the edges passing through any vertex
of the $K_{l+1}$ subgraph will form a linearly independent set of
$l+1$ elements. Again we get a contradiction with the fact that
$r(\mathcal A_X)=l$. So $G_X$ is a $K_{l+1}$.

\vskip .2in

With this, $G_X$ has exactly $l+1$ $K_l$ subgraphs. These
subgraphs will give us the 'special' elements of $D_{l-1}(\mathcal
A_X)$: $r_{2,\ldots,l+1},r_{1,3,\ldots,l+1}, \ldots,
r_{1,2,\ldots,l}$. The only relation on these elements is exactly
the 'special' element in $R_l(\mathcal A_X)$:
$$r_{2,\ldots,l+1}-r_{1,3,\ldots,l+1}+\cdots +
(-1)^lr_{1,2,\ldots,l}.$$ We denote this element with
$r_{1,2,\ldots,l,l+1}$ and he is forming the basis for
$R_l(\mathcal A_X)$.

\vskip .2in

For the converse, it is obvious that if $G_X$ is a $K_{l+1}$, then
$R_l(\mathcal A_X)\neq 0$ and even more, $\dim R_l(\mathcal
A_X)=1$.
\end{proof}

With this lemma we can identify easily the two complexes. The way
we pick the special basis elements will give us the same matrices
for the differentials of the two complexes and, hence, with Lemma
2.5., we have proved the following proposition:

\begin{prop} Let $G$ be a connected graph. $\mathcal A_G$ is $k-$formal if and
only if $H_i(\Delta)=0$ for every $i=1,\ldots,k-1$.
\end{prop}

Note that from this proposition we get that in the graphic
arrangement case, $k-$formality depends only on combinatorics,
contrary to the case of lines arrangements (see Yuzvinsky's
example).

\begin{exm} We conclude with an easy example of a formal graphic
arrangement which is not 3-formal. Consider the graph $G$ in the
figure below:
\begin{figure}[h]
\begin{center}
\includegraphics{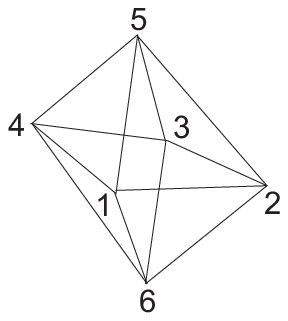}
\end{center}

\end{figure}

The associated flag complex $\Delta$ is the boundary complex of an
octahedron on the same vertices and edges. The associated chain
complex of $\Delta$ is:
$$0 \longrightarrow \mathbb K^8 \stackrel{f_2} {\longrightarrow}
\mathbb K^{12} \stackrel{f_1}{\longrightarrow} \mathbb K^6
\longrightarrow 0,$$ where, if we order the basis
lexicographically we have:

$$f_1=\left[
\begin{array}{cccccccccccc}
1&1&1&1&0&0&0&0&0&0&0&0 \\
-1&0&0&0&1&1&1&0&0&0&0&0 \\
0&0&0&0&-1&0&0&1&1&1&0&0 \\
0&-1&0&0&0&0&0&-1&0&0&1&1 \\
0&0&-1&0&0&-1&0&0&-1&0&-1&0 \\
0&0&0&-1&0&0&-1&0&0&-1&0&-1
\end{array}
\right]$$ and

$$f_2=\left[
\begin{array}{cccccccc}
1&1&0&0&0&0&0&0 \\
0&0&1&1&0&0&0&0 \\
-1&0&-1&0&0&0&0&0 \\
0&-1&0&-1&0&0&0&0 \\
0&0&0&0&1&1&0&0 \\
1&0&0&0&-1&0&0&0 \\
0&1&0&0&0&-1&0&0 \\
0&0&0&0&0&0&1&1 \\
0&0&0&0&1&0&-1&0 \\
0&0&0&0&0&1&0&-1 \\
0&0&1&0&0&0&1&0 \\
0&0&0&1&0&0&0&1
\end{array}
\right].$$

Since $G$ is connected, then $\dim H_0(\Delta)=1$. So
$rk(f_1)=6-1=5$. Therefore, $\dim ker(f_1)=12-5=7$.

Every 4-cycle in $G$ is a linear combination of 3-cycles. So
$\mathcal A_G$ is formal (2-formal). By the proposition above,
$\dim H_1(\Delta)=0$ and with this we get $rk(f_2)=7$. Therefore,
$\dim ker(f_2)=8-7=1$. So we get $\dim H_2(\Delta)=1$. Hence
$\mathcal A_G$ is not 3-formal.
\end{exm}

\vskip 0.3in

\noindent{\bf Acknowledgment} I thank Hal Schenck for many useful
conversations and grant support.

\renewcommand{\baselinestretch}{1.0}
\small\normalsize 

\bibliographystyle{amsalpha}

\end{document}